\documentclass{amsart}

\usepackage{amssymb}
\usepackage{amsthm}
\usepackage{amsmath}
\usepackage[hidelinks]{hyperref}

\newcommand{\numberseries}{\bfseries}   

\newlength{\thmtopspace}                
\newlength{\thmbotspace}                
\newlength{\thmheadspace}               
\newlength{\thmindent}                  

\renewcommand{\subparagraph}{\vspace*{\thmbotspace}}

\newtheoremstyle{bfupright head,slanted body}
                {\thmtopspace}{\thmbotspace}
                {\slshape}{\thmindent}{\bfseries}{.}{\thmheadspace}
                {{\numberseries \thmnumber{#2\;}}\thmnote{#3}}

\newtheoremstyle{bfupright head,upright body}
                {\thmtopspace}{\thmbotspace}
                {\upshape}{\thmindent}{\bfseries}{.}{\thmheadspace}
                {{\numberseries \thmnumber{#2\;}}\thmnote{#3}}

\newtheoremstyle{bfit head,upright body}
                {\thmtopspace}{\thmbotspace}
                {\upshape}{\thmindent}{\upshape}{.}{\thmheadspace}
                {{\numberseries\thmnumber{#2\;}}
                {\bfseries\itshape\thmnote{\negthickspace#3}}}

\newtheoremstyle{it head,upright body}
                {\thmtopspace}{\thmbotspace}
                {\upshape}{\thmindent}{\upshape}{.}{\thmheadspace}
                {{\numberseries\thmnumber{#2\;}}
                {\itshape\thmnote{\negthickspace#3}}}

\newtheoremstyle{fixed bf head,slanted body}
                {\thmtopspace}{\thmbotspace}{\slshape}
                {\thmindent}{\bfseries}{.}{\thmheadspace}
                {{\numberseries \thmnumber{#2\;}}\thmname{#1}\thmnote{ (#3)}}

\newtheoremstyle{fixed bf head,upright body}
                {\thmtopspace}{\thmbotspace}{\upshape}
                {\thmindent}{\bfseries}{.}{\thmheadspace}
                {{\numberseries \thmnumber{#2\;}}\thmname{#1}\thmnote{ (#3)}}

\newtheoremstyle{fixed bfit head,upright body}
                {\thmtopspace}{\thmbotspace}{\upshape}
                {\thmindent}{\bfseries\itshape}{.}{\thmheadspace}
                {{\numberseries \thmnumber{#2\;}}\thmname{#1}\thmnote{ (#3)}}

\newtheoremstyle{sc head,small body}
                {\thmtopspace}{\thmbotspace}
                {\small\upshape}{\thmindent}{\scshape}{.}{\thmheadspace}
                {\thmname{#1}}

\newtheoremstyle{numbered paragraph}
                {\thmtopspace}{\thmbotspace}{\upshape}
                {\thmindent}{\upshape}{}{\thmheadspace}
                {{\numberseries \thmnumber{#2.}}}

\newtheoremstyle{unnumbered paragraph}
                {\thmtopspace}{\thmbotspace}{\upshape}
                {\parindent}{\upshape}{}{0pt}

\theoremstyle{bfupright head,slanted body}
\newtheorem{res}{}[section]             \newtheorem*{res*}{}

\theoremstyle{fixed bf head,slanted body}
\newtheorem{thm}[res]{Theorem}          \newtheorem*{thm*}{Theorem}
\newtheorem{prp}[res]{Proposition}      \newtheorem*{prp*}{Proposition}
\newtheorem{cor}[res]{Corollary}        \newtheorem*{cor*}{Corollary}
\newtheorem{lem}[res]{Lemma}            \newtheorem*{lem*}{Lemma}

\theoremstyle{fixed bf head,upright body}
       \newtheorem*{dfn*}{Definition}
     \newtheorem*{con*}{Construction}
\newtheorem{rmk}[res]{Remark}           \newtheorem*{rmk*}{Remark}
\newtheorem{exa}[res]{Example}          \newtheorem*{exa*}{Example}

\theoremstyle{numbered paragraph}
\newtheorem{ipg}[res]{}

\newlength{\thmlistleft}        
\newlength{\thmlistright}       
\newlength{\thmlistpartopsep}   
\newlength{\thmlisttopsep}      
\newlength{\thmlistparsep}      
\newlength{\thmlistitemsep}     

\newcounter{eqc} 
\newenvironment{eqc}{\begin{list}{\upshape (\textit{\roman{eqc}})}%
    {\usecounter{eqc}%
      \setlength{\leftmargin}{\thmlistleft}%
      \setlength{\labelwidth}{\thmlistleft}%
      \setlength{\rightmargin}{\thmlistright}%
      \setlength{\partopsep}{\thmlistpartopsep}%
      \setlength{\topsep}{\thmlisttopsep}%
      \setlength{\parsep}{\thmlistparsep}%
      \setlength{\itemsep}{\thmlistitemsep}}}%
  {\end{list}}%

\newcommand{\eqclbl}[1]{{\upshape(\textit{#1})}}

\newcounter{prt}
\newenvironment{prt}{\begin{list}{\upshape (\alph{prt})}%
    {\usecounter{prt}%
      \setlength{\leftmargin}{\thmlistleft}%
      \setlength{\labelwidth}{\thmlistleft}%
      \setlength{\rightmargin}{\thmlistright}%
      \setlength{\partopsep}{\thmlistpartopsep}%
      \setlength{\topsep}{\thmlisttopsep}%
      \setlength{\parsep}{\thmlistparsep}%
      \setlength{\itemsep}{\thmlistitemsep}}}%
  {\end{list}}%

\newcommand{\prtlbl}[1]{{\upshape(#1)}}

\newcounter{rqm}
\newenvironment{rqm}{\begin{list}{\upshape (\arabic{rqm})}%
    {\usecounter{rqm}%
      \setlength{\leftmargin}{\thmlistleft}%
      \setlength{\labelwidth}{\thmlistleft}%
      \setlength{\rightmargin}{\thmlistright}%
      \setlength{\partopsep}{\thmlistpartopsep}%
      \setlength{\topsep}{\thmlisttopsep}%
      \setlength{\parsep}{\thmlistparsep}%
      \setlength{\itemsep}{\thmlistitemsep}}}%
  {\end{list}}%
\newenvironment{prf*}[1][Proof]{%
  \begin{proof}[\bf #1]
    \setcounter{equation}{0}
    }
  {\end{proof}
}

  \newcommand{\proofoftag}[2][:]{(#2)#1}

  \newcommand{\proofofimp}[3][:]{\mbox{\eqclbl{#2}$\!\implies\!$\eqclbl{#3}#1}}

  \newcommand{\pgref}[1]{\ref{#1}}
  \renewcommand{\eqref}[1]{(\pgref{eq:#1})}
  \newcommand{\thmref}[2][Theorem~]{#1\pgref{thm:#2}}
  \newcommand{\corref}[2][Corollary~]{#1\pgref{cor:#2}}
  \newcommand{\prpref}[2][Proposition~]{#1\pgref{prp:#2}}
  \newcommand{\lemref}[2][Lemma~]{#1\pgref{lem:#2}}
  \newcommand{\exaref}[2][Example~]{#1\pgref{exa:#2}}
  \newcommand{\rmkref}[2][Remark~]{#1\pgref{rmk:#2}}
  \newcommand{\secref}[2][Section~]{#1\ref{sec:#2}}
  \newcommand{\thmcite}[2][?]{\cite[Thm.~#1]{#2}}
  \newcommand{\corcite}[2][?]{\cite[Cor.~#1]{#2}}
  \newcommand{\prpcite}[2][?]{\cite[Prop.~#1]{#2}}
  \newcommand{\lemcite}[2][?]{\cite[Lem.~#1]{#2}}
  \newcommand{\seccite}[2][?]{\cite[Sect.~#1]{#2}}
  \newcommand{\exacite}[2][?]{\cite[Example~#1]{#2}}

\def\urltilda{\kern -.15em\lower .7ex\hbox{\~{}}\kern .04em} 

\newcommand{\ZZ}{\mathbb{Z}}
\newcommand{\dis}{\:\is\:}
\newcommand{\fm}{\mathfrak{m}}
\newcommand{\is}{\cong}
\newcommand{\qis}{\simeq}
\newcommand{\onto}{\twoheadrightarrow}
\newcommand{\into}{\hookrightarrow}
\newcommand{\lra}{\longrightarrow}
\newcommand{\xra}[2][]{\xrightarrow[#1]{\;#2\;}}
\newcommand{\qra}{\xra{\qis}}
\newcommand{\QQ}{\mathbb{Q}}
\newcommand{\QZ}{\QQ/\ZZ}
\newcommand{\Rop}{R^\circ}
\newcommand{\Cy}[2][]{\operatorname{Z}_{#1}(#2)}
\renewcommand{\H}[2][]{\operatorname{H}_{#1}(#2)}
\newcommand{\Shift}[2][]{\mathsf{\Sigma}^{#1}{#2}}

\newcommand{\Hom}[3][R]{\operatorname{Hom}_{#1}(#2,#3)}
\newcommand{\Ext}[4][R]{\operatorname{Ext}_{#1}^{#2}(#3,#4)}
\newcommand{\tp}[3][R]{\nobreak{#2\otimes_{#1}#3}}
\newcommand{\Tor}[4][R]{\operatorname{Tor}^{#1}_{#2}(#3,#4)}
\providecommand{\Mod}[1][R]{\Cat{#1}{M}}

\newcommand{\maintheorem}{Theorem~0}
\newcommand{\Prj}[1][R]{\mathsf{Prj}(#1)}
\newcommand{\Inj}[1][R]{\mathsf{Inj}(#1)}
\newcommand{\FPPrj}[1][\Rop]{\mathsf{FpPrj}(#1)}
\newcommand{\FPInj}[1][\Rop]{\mathsf{FpInj}(#1)}
\newcommand{\Flat}[1][R]{\mathsf{Flat}(#1)}
\newcommand{\Cot}[1][R]{\mathsf{Cot}(#1)}
\renewcommand{\Mod}[1][R]{\mathsf{Mod}(#1)}
\newcommand{\semi}[1]{\mathsf{semi\textnormal{-}#1}}
\newcommand{\ac}[1]{\mathsf{#1\textnormal{-}ac}}
 
\newcommand{\tagform}[1]{\textnormal{\small(#1)}}
\newcommand{\tagrange}[2]{\textnormal{\small(#1)--(#2)}}

\newcommand{\tagPf}{\tagform{P0}}
\newcommand{\tagP}{\tagform{P1}}
\newcommand{\tagPc}{\tagform{P2}}
\newcommand{\tagPa}{\tagform{P3}}
\newcommand{\tagsP}{\tagrange{P0}{P3}}

\newcommand{\tagF}{\tagform{F1}}
\newcommand{\tagFp}{\tagform{F2}}
\newcommand{\tagFa}{\tagform{F3}}
\newcommand{\tagsF}{\tagrange{F1}{F3}}

\newcommand{\tagFC}{\tagform{FC1}}
\newcommand{\tagFCc}{\tagform{FC2}}
\newcommand{\tagFCa}{\tagform{FC3}}
\newcommand{\tagsFC}{\tagrange{FC1}{FC3}}

\newcommand{\tagI}{\tagform{I1}}
\newcommand{\tagIc}{\tagform{I2}}
\newcommand{\tagIa}{\tagform{I3}}
\newcommand{\tagsI}{\tagrange{I1}{I3}}

\newcommand{\tagFPI}{\tagform{fpI1}}
\newcommand{\tagFPIp}{\tagform{fpI2}}
\newcommand{\tagFPIa}{\tagform{fpI3}}
\newcommand{\tagsFPI}{\tagrange{fpI1}{fpI3}}

\newcommand{\tagFPPI}{\tagform{fpPI1}}
\newcommand{\tagFPPIc}{\tagform{fpPI2}}
\newcommand{\tagFPPIa}{\tagform{fpPI3}}
\newcommand{\tagsFPPI}{\tagrange{fpPI1}{fpPI3}}

\numberwithin{equation}{res}

\title{Acyclic complexes and regular rings}

\author[L.W.\ Christensen]{Lars Winther Christensen} %
\address{L.W.C. \ Texas Tech University, Lubbock, TX 79409, U.S.A.}
\email{lars.w.christensen@ttu.edu}
\urladdr{http://www.math.ttu.edu/\urltilda lchriste}

\author[S.\ Estrada]{Sergio Estrada} %
\address{S.E. \ Universidad de Murcia, Murcia 30100, Spain}
\email{sestrada@um.es} %
\urladdr{https://webs.um.es/sestrada/}

\author[P.\ Thompson]{Peder Thompson}

\address{P.T. \ M\"{a}lardalen University, V\"{a}ster{\aa}s 72123,
  Sweden}

\email{peder.thompson@mdu.se}

\urladdr{https://sites.google.com/view/pederthompson}

\thanks{L.W.C.\ was partly supported by Simons Foundation
  collaboration grant 962956.  S.E. was partly supported by grant
  22004/PI/22 funded by Fundaci\'on S\'eneca, Agencia de Ciencia y
  Tecnolog\'ia de la Regi\'on de Murcia and by grants
  PID2020-113206GB-I00 and PID2024-155576NB-I00 funded by MICIU/ AEI/10.13039/501100011033.}

\date{21 January 2026}



\begin{document}

\begin{abstract}
  A 2009 paper by Iacob and Iyengar characterizes noetherian regular
  rings in terms of properties of complexes of projective modules,
  flat modules, and injective modules. We show that the relevant
  properties of such complexes are equivalent without reference to
  regularity of the ring and that they characterize coherent regular
  rings and von Neumann regular rings.
\end{abstract}

\maketitle


\section*{Introduction}

\noindent
In this short paper $R$ denotes an associative unital ring. An
$R$-module is a left $R$-module, and right $R$-modules are considered
modules over the opposite ring $\Rop$.

Following Bertin \cite{JBr71} and Glaz \cite{ccr} we say that $R$ is
left/right \emph{regular} if every finitely generated left/right ideal
in $R$ has finite projective dimension.  We note that this definition
is broader than the one used by Iacob and Iyengar \cite{AIcSBI09},
which implicitly includes the assumption that $R$ is left/right
noetherian.

As a corollary to our first main result, \thmref{P}, we obtain the
following:%

\begin{res*}[\maintheorem]
  Let $R$ be right coherent. The following conditions are equivalent.
  \setlength{\thmlistleft}{3em} %
  \begin{eqc}
  \item $R$ is right regular.
  \item Every complex of finitely generated free $R$-modules is
    semi-projective.
  \item Every complex of projective $R$-modules is semi-projective.
  \item Every complex of injective $\Rop$-modules is semi-injective.
  \item Every complex of flat $R$-modules is semi-flat.
  \item Every acyclic complex of projective $R$-modules is
    contractible.
  \item Every acyclic complex of injective $\Rop$-modules is
    contractible.
  \item Every acyclic complex of flat $R$-modules is pure acyclic.
  \end{eqc}
  \setlength{\thmlistleft}{2.5em} %
\end{res*}

\noindent In fact, \thmref{P} shows that any right regular ring
satisfies the conditions above.

The equivalence of the conditions in \maintheorem\ was proved for
commutative noetherian rings by Christensen, Foxby, and Holm
\cite{dcmca}. Beyond that realm it applies, for example, to von
Neumann regular rings; see \exaref{vnr}.

The equivalence of some of the conditions in \maintheorem\ was proved
already in \cite{AIcSBI09}, see \rmkref{1}. Recently the equivalence
of \eqclbl{i} and \eqclbl{iii}--\eqclbl{viii} was proved by Gillespie
and Iacob \cite{JGlAIc}; see \rmkref{2}. Our proofs do not rely on
these works.

In \secref[Sections~]{2} and \secref[]{3} we characterize right
regularity of a right coherent ring $R$ in terms of properties of
complexes of flat-cotorsion $R$-modules and complexes of fp-injective
$\Rop$-modules. One upshot is that the equivalence of conditions
\eqclbl{ii}--\eqclbl{viii} in \maintheorem\ can be established without
invoking \eqclbl{i}, the right regularity of the ring; this is part of
our second main result, \thmref{FP}.

\section*{Acknowledgment}

\noindent
This work was done during a joint visit by L.W.C.\ and S.E.\ to
M\"{a}lardalen University, and the institution's hospitality is
acknowledged with gratitude. Further, we thank the referee for recommendations that
improved the presentation.

\section{Complexes of projective, injective, and flat modules}
\label{sec:1}

\noindent
In this paper complexes, i.e.\ complexes of modules, are indexed
homologically:
\begin{equation*}
  X = \cdots \lra X_{n+1} \lra X_n \lra X_{n-1} \lra \cdots \:.
\end{equation*}
For $s \in \ZZ$ the complex $\Shift[s]{X}$ has
$(\Shift[s]{X})_n = X_{n-s}$ and differential
$\partial^{\Shift[s]{X}} = (-1)^s\partial^X$.  For convenience we give references to
\cite{dcmca}, which the reader may also consult for any unexplained
notation or terminology.

We start by recalling some facts about complexes of
projective and flat modules.

\begin{ipg}
  \label{proj}
  Recall from \prpcite[5.2.10]{dcmca} that a complex $P$ of projective
  $R$-modules is semi-projective if and only if $\Hom{P}{A}$ is
  acyclic for every acyclic $R$-complex~$A$.

  For every $R$-complex $X$ there is an exact sequence of
  $R$-complexes
  \begin{equation}
    \label{eq:proj}
    0 \lra A \lra P \lra X \lra 0    
  \end{equation}
  where $P$ is semi-projective and $A$ is acyclic.
  \begin{prt}
  \item For an acyclic complex $P$ of projective $R$-modules the next
    conditions are equivalent. (See \prpcite[4.3.29 and
    Cor.~5.5.26]{dcmca}.)
    \begin{eqc}
    \item $P$ is semi-projective.
    \item $P$ is contractible.
    \item $P$ is pure acyclic.
    \end{eqc}
  \item If in an exact sequence $0 \to P' \to P \to P'' \to 0$ of
    complexes of projective $R$-modules two of the complexes are
    semi-projective, then so is the third.
  \end{prt}
\end{ipg}

In the context of the derived category over $R$ the semi-projective
$R$-complexes recalled above are used to define derived Hom and tensor
product functors. Indeed, for every $R$-complex $X$ there is a
quasi-isomorphism $P \qra X$ with $P$ semi-projective; this is known
as a semi-projective resolution and encompasses the classic notion of
a projective resolution. The derived tensor product functor can also
be computed using the semi-flat complexes recalled below, and the
(contravariant) derived Hom functor can be defined using the
semi-injective complexes recalled in \pgref{inj}.

\enlargethispage*{2\baselineskip}
\begin{ipg}
  \label{flat}
  Recall from \prpcite[5.4.9]{dcmca} that a complex $F$ of flat
  $R$-modules is semi-flat if and only if $\tp{A}{F}$ is acyclic for
  every acyclic $\Rop$-complex $A$.
  \begin{prt}
  \item An acyclic complex of flat $R$-modules is pure acyclic if and
    only if it is semi-flat. (See \thmcite[5.5.22]{dcmca}.)
  \item If in an exact sequence $0 \to F' \to F \to F'' \to 0$ of
    complexes of flat $R$-modules two of the complexes are semi-flat,
    then so is the third.
  \item A complex of projective $R$-modules is semi-flat if and only
    if it is semi-projective. (See \corcite[5.4.10 and
    Thm.~5.5.27]{dcmca}.)
  \end{prt}
\end{ipg}

\begin{prp}
  \label{prp:P}
  The following conditions are equivalent.
  \begin{rqm}
  \item[\tagPf] Every complex of finitely generated free $R$-modules
    is semi-projective.
  \item[\tagP] Every complex of projective $R$-modules is
    semi-projective.
  \item[\tagPc] Every acyclic complex of projective $R$-modules is
    contractible.
  \item[\tagPa] Every acyclic complex of projective $R$-modules is
    semi-projective.

  \item[\tagF] Every complex of flat $R$-modules is semi-flat.
  \item[\tagFp] Every acyclic complex of flat $R$-modules is pure
    acyclic.
  \item[\tagFa] Every acyclic complex of flat $R$-modules is
    semi-flat.
  \end{rqm}
\end{prp}

\begin{prf*}
  First we argue that conditions \tagsP\ are equivalent.  In view of
  \pgref{proj}(a) and the fact that free modules are projective, the
  following implications are clear:
  \begin{equation*}
    \tagPf \impliedby \tagP \implies \tagPc \implies \tagPa \,,
  \end{equation*}
  which leaves two implications to prove.

  \tagPa$\implies$\tagP: Let $X$ be a complex of projective
  $R$-modules and consider the exact sequence
  $0 \to A \to P \to X \to 0$ from \eqref{proj}. It follows that $A$
  is a complex of projective $R$-modules and, thus, semi-projective,
  so $X$ is semi-projective by~\pgref{proj}(b).

  \tagPf$\implies$\tagP: Let $\mathcal{S}$ be a set of representatives
  for the isomorphism classes of bounded above complexes of finitely
  generated free $R$-modules; it generates a cotorsion pair
  $({}^\perp(\mathcal{S}^\perp),{S}^\perp)$ in the category of
  $R$-complexes. By work of Bravo, Gillespie, and Hovey
  \thmcite[A.3]{BGH} the left-hand class ${}^\perp(\mathcal{S}^\perp)$
  consists of all complexes of projective $R$-modules. Thus, if every
  complex in $\mathcal{S}$ is semi-projective, then every acyclic
  $R$-complex belongs to $\mathcal{S}^\perp$, whence every complex of
  projective $R$-modules is semi-projective.

  Condition \tagF\ evidently implies \tagFa. By \pgref{flat}(a)
  conditions \tagFp\ and \tagFa\ are equivalent, and by
  \pgref{proj}(a) they imply \tagPc. This leaves one implication to
  prove.

  \tagP$\implies$\tagF: Let $F$ be a complex of flat $R$-modules; by
  work of Neeman \thmcite[8.6]{ANm08} there is an exact sequence
  \begin{equation*}
    0 \lra A \lra P \lra F \lra 0
  \end{equation*}
  where $P$ is a complex of projective $R$-modules and $A$ is a pure
  acyclic complex of flat $R$-modules. By assumption $P$ is
  semi-projective and hence semi-flat, see \pgref{flat}(c). The
  complex $A$ is semi-flat by \pgref{flat}(a), so $F$ is semi-flat by
  \pgref{flat}(b).
\end{prf*}

\begin{ipg}
  \label{inj}
  Recall from \prpcite[5.3.16]{dcmca} that a complex $I$ of injective
  $\Rop$-modules is semi-injective if and only if $\Hom[\Rop]{A}{I}$
  is acyclic for every acyclic $\Rop$-complex $A$.

  For every $\Rop$-complex $X$ there is an exact sequence of
  $\Rop$-complexes
  \begin{equation}
    \label{eq:inj}
    0 \lra X \lra I \lra A \lra 0    
  \end{equation}
  where $I$ is semi-injective and $A$ is acyclic.
  \begin{prt}
  \item For an acyclic complex $I$ of injective $\Rop$-modules the
    following conditions are equivalent. (See \prpcite[4.3.29]{dcmca}
    and Bazzoni, Cort\'{e}s-Izurdiaga, and Estrada \cite[Props.~2.4(1)
    and 4.8(1)]{BCE-20}.)
    \begin{eqc}
    \item $I$ is semi-injective.
    \item $I$ is contractible.
    \item $I$ is pure acyclic.
    \end{eqc}
  \item If in an exact sequence $0 \to I' \to I \to I'' \to 0$ of
    complexes of injective $\Rop$-modules, two of the complexes are
    semi-injective, then so is the third.
  \end{prt}
\end{ipg}
\begin{prp}
  \label{prp:I}
  The following conditions are equivalent.
  \begin{rqm}
  \item[\tagI] Every complex of injective $\Rop$-modules is
    semi-injective.
  \item[\tagIc] Every acyclic complex of injective $\Rop$-modules is
    contractible.
  \item[\tagIa] Every acyclic complex of injective $\Rop$-modules is
    semi-injective.
  \end{rqm}
\end{prp}

\begin{prf*}
  By \pgref{inj}(a) conditions \tagIc\ and \tagIa\ are equivalent, and
  \tagI\ evidently implies \tagIa. For the converse let $X$ be a
  complex of injective $\Rop$-modules and consider the exact sequence
  $0 \to X \to I \to A \to 0$ from \eqref{inj}. It follows that $A$ is
  a complex of injective $\Rop$-modules and, therefore,
  semi-injective, so $X$ is semi-injective by~\pgref{inj}(b).
\end{prf*}

\begin{rmk}
  \label{rmk:0}
  The conditions \tagsI\ imply the conditions from \prpref{P}.  To see
  this it suffices to verify that \tagIa\ implies \tagFp: Let $F$ be
  an acyclic complex of flat $R$-modules. The character complex
  $\Hom[\ZZ]{F}{\QZ}$ is an acyclic complex of injective
  $\Rop$-modules and, therefore, contractible, whence $F$ is pure
  acyclic.
\end{rmk}

\begin{rmk}
  \label{rmk:1}
  The equivalence of some of the conditions above are known from the
  literature. Conditions \tagP, \tagPc, \tagF, and \tagFp\ are
  equivalent by \cite[Props.~3.1, 3.3, and 3.4]{AIcSBI09}. Conditions
  \tagI\ and \tagIc\ are equivalent by \prpcite[2.1]{AIcSBI09}.

  Under the assumption that $R$ is commutative and noetherian, the
  equivalence of \tagPf\ and \tagP\ was proved in
  \thmcite[20.2.12]{dcmca}. Tereshkin \cite{DTr} has informed us that
  he knows of the equivalence for any ring, as proved in \prpref{P},
  from private communication with Positselski. The argument he implied
  is different from the one we give here.
\end{rmk}

\section{Regular rings}

\begin{thm}
  \label{thm:P}
  Consider the following conditions on
  $R$. \setlength{\thmlistleft}{3em} %
  \begin{eqc}
  \item $R$ is right regular.
  \item Every finitely generated right ideal in $R$ has finite flat
    dimension.
  \item One/all of conditions \tagsI\ from {\rm\prpref[]{I}} hold.
  \item One/all of conditions \tagsP\ and \tagsF\ from
    {\rm\prpref[]{P}} hold.
  \item Every $\Rop$-module with a degreewise finitely generated
    projective resolution has finite projective dimension.
  \end{eqc} \setlength{\thmlistleft}{2.5em} %
  The following implications hold,
  \begin{center}
    \eqclbl{i} $\implies$ \eqclbl{ii} $\implies$ \eqclbl{iii}
    $\implies$ \eqclbl{$iv$} $\implies$ \eqclbl{v}\,,
  \end{center}
  and if $R$ is right coherent then \eqclbl{v} implies \eqclbl{i},
  i.e.\ the conditions are equivalent.
\end{thm}

\begin{prf*}
  Evidently \eqclbl{i} implies \eqclbl{ii}, and \eqclbl{iii} implies
  \eqclbl{iv} by \rmkref{0} .

  \proofofimp{ii}{iii} By \prpref{I} it suffices to show that $R$
  satisfies \tagIc. Let $I$ be an acyclic complex of injective
  $\Rop$-modules and set $J = \prod_{n\in\ZZ}\Shift[n]{I}$. It is also
  an acyclic complex of injective $\Rop$-modules, and one has
  $\Cy[n]{J} \is \Cy[n-1]{J}$ for every $n\in\ZZ$, so it follows from
  \prpcite[4.8(3)]{BCE-20} that the cycle modules $\Cy[n]{J}$ are
  injective. Each cycle module $\Cy[n]{I}$ is a direct summand of
  $\Cy[n]{J}$ and hence injective, so $I$ is contractible.
  
  \proofofimp{iv}{v} By \prpref{P} it suffices to show that \tagsP\
  imply that every $\Rop$-module with a degreewise finitely generated
  projective resolution has finite projective dimension. Let
  $L \qra M$ be such a module and resolution. Further, let $E$ be a
  faithfully injective $R$-module and $P \qra E$ a projective
  resolution. The complex $L^* = \Hom[\Rop]{L}{R}$ of finitely
  generated projective $R$-modules is by \tagPf\ semi-projective; this
  explains the first and third isomorphisms in the computation
  below. Further, $L^*$ is by \prpcite[7.12]{ANm08} a compact object
  in the homotopy category of projective $R$-modules; this explains
  the second isomorphism.
  \begin{align*}
    \textstyle \H[0]{\Hom{L^*}{\coprod_{n \in \ZZ}\Shift[n]{E}}} 
    & \textstyle \dis \H[0]{\Hom{L^*}{\coprod_{n \in \ZZ}\Shift[n]{P}}} \\
    & \textstyle \dis \coprod_{n \in \ZZ}\H[0]{\Hom{L^*}{\Shift[n]{P}}} \\
    & \textstyle \dis \coprod_{n \in \ZZ}\H[0]{\Hom{L^*}{\Shift[n]{E}}} \\    
    & \textstyle \dis \coprod_{n \in \ZZ}{\H[n]{\Hom{L^*}{E}}}  \\
    & \textstyle \dis \coprod_{n \in \ZZ}{\Hom{\H[n]{L^*}}{E}}  \,.      
  \end{align*}
  At the same time, the canonical embedding
  $\coprod_{n \in \ZZ}\Shift[n]{E} \to \prod_{n \in \ZZ}\Shift[n]{E}$
  is an isomorphism, so one has
  \begin{align*}
    \textstyle \H[0]{\Hom{L^*}{\coprod_{n \in \ZZ}\Shift[n]{E}}} 
    & \textstyle \dis \H[0]{\Hom{L^*}{\prod_{n \in \ZZ}\Shift[n]{E}}} \\
    & \dis \textstyle  \prod_{n \in \ZZ}\H[0]{\Hom{L^*}{\Shift[n]{E}}} \\
    & \dis \textstyle \prod_{n \in \ZZ} \H[n]{\Hom{L^*}{E}} \\
    & \dis \textstyle \prod_{n \in \ZZ} \Hom{\H[n]{L^*}}{E} \,.      
  \end{align*}
  As the relevant Hom, homology, and shift functors preserve
  (co)products in the sense of \cite[3.1.8 and 3.1.20]{dcmca}, it
  follows that the canonical embedding
  \begin{equation*}
    \textstyle \coprod_{n\in\ZZ} \Hom{\H[n]{L^*}}{E} \lra \prod_{n\in\ZZ} \Hom{\H[n]{L^*}}{E}
  \end{equation*}
  is an isomorphism. Thus, $\Hom{\H[n]{L^*}}{E}$ and, therefore,
  $\H[n]{L^*}$ is non-zero for only finitely many $n \in \ZZ$. Thus,
  for $n \ll 0$ the complex $0 \to \Cy[n]{L^*} \to (L^*)_n \to \cdots$
  is acyclic, so splicing it with a projective resolution of
  $\Cy[n]{L^*}$ yields an acyclic complex of projective
  $R$-modules. By \tagPa\ it is contractible, so for $n \ll 0$ the
  $R$-module $\Cy[n]{L^*}$ is projective. As one has
  $L \is \Hom{L^*}{R}$ the $\Rop$-module $\Cy[n]{L}$ is projective for
  $n \gg 0$, so $M$ has finite projective dimension.

  \proofofimp{v}{i} Assume that $R$ is right coherent and let
  $\mathfrak{a} \subset R$ be a finitely generated right ideal. The
  quotient $R/\mathfrak{a}$ is a finitely presented $\Rop$-module, and
  every such module has a degreewise finitely generated projective
  resolution. It follows that $R/\mathfrak{a}$ and, therefore,
  $\mathfrak{a}$ has finite projective dimension.
\end{prf*}

Without the coherence assumption the last implication proved above may
fail.

\begin{exa}
  Let $k$ be a field. The local ring
  $R = k[x_1,x_2,\ldots]/(x_1,x_2,\ldots)^2$ with maximal ideal
  $\fm = (x_1,x_2,\ldots)$ is not coherent; indeed the kernel $\fm$ of
  the canonical map $R \onto (x_1)$ is not finitely generated. By
  \prpcite[2.5]{BGH} the only $R$-modules that admit a degreewise
  finitely generated projective resolution are the finitely generated
  free $R$-modules. However, $R$ is not a regular ring: The proof of
  \prpcite[2.5]{BGH} shows that every non-free finitely presented
  $R$-module has projective dimension at least~$2$. Yet, the existence
  of any $R$-module of finite projective dimension at least $1$
  implies the existence of an injective homomorphism
  $\partial\colon P \into Q$ of free $R$-modules. As $R$ is perfect,
  one can assume that the image of $\partial$ is contained in $\fm Q$,
  see \cite[Thms.~B.55 and B.60]{dcmca}, which forces the existence of
  a non-free finitely presented $R$-module of projective
  dimension~$1$; a contradiction. Thus, every non-projective
  $R$-module has infinite projective dimension.
\end{exa}

The example above suggests that coherence is the ``minimal'' condition
on $R$ that makes all five conditions in \thmref{P} equivalent. In the
proof of \thmref{P}, the argument for the implication
\proofofimp[]{iv}{v} relies on Neeman's \prpcite[7.12]{ANm08}, which
is also used in \cite[3.5]{AIcSBI09}. Let us, therefore, record that
this implication, under the coherence assumption, has a more
elementary proof.

\begin{rmk}
  Let $R$ be right coherent and assume that it satisfies
  \thmref{P}\eqclbl{iv}. Let $M$ be a finitely presented
  $\Rop$-module, and $N$ an $R$-module with a semi-flat resolution
  $F \qra N$. The canonical embedding
  $\Phi\colon\coprod_{n \in \ZZ} \Shift[n]{F} \to \prod_{n \in \ZZ}
  \Shift[n]{F}$ is a quasi-isomorphism since the homology of either
  complex equals $N$ in each degree. The mapping cone of $\Phi$ is by
  $\tagFp$ pure acyclic, so $\tp{M}{\Phi}$ is a quasi-isomorphism as
  well; this explains the third isomorphism below.
  \begin{align*}
    \textstyle\coprod_{n \in \ZZ}\Tor{n}{M}{N} 
    & \dis \textstyle\coprod_{n \in \ZZ}\H[0]{\tp{M}{\Shift[n]{F}}} \\%
    & \dis \textstyle\H[0]{\tp{M}{\coprod_{n \in \ZZ}\Shift[n]{F}}} \\%
    & \dis \textstyle\H[0]{\tp{M}{\prod_{n \in \ZZ}\Shift[n]{F}}} \\%
    & \dis \textstyle\prod_{n \in \ZZ}\H[0]{\tp{M}{\Shift[n]{F}}} 
    \dis \textstyle\prod_{n \in \ZZ}\Tor{n}{M}{N} 
  \end{align*}
  As in the proof of \thmref{P} it follows that $\Tor{n}{M}{N}$ is
  nonzero for only finitely many $n$ whence $M$, being finitely
  presented, has finite projective dimension. Thus $R$ satisfies
  \thmref{P}\eqclbl{v}.
\end{rmk}

The proof of \thmref{P} suggests that the equivalence, for a right
coherent ring, of all ten conditions from \prpref[Propositions~]{P}
and \prpref[]{I} ``factors through'' the right regularity property of
the ring, but their equivalence can, in fact, be established without
reference to this property; see \rmkref{3} and \thmref{FP}.

As an application of \thmref{P} we offer a short proof of a result
already proved by Glaz \thmcite[6.2.5]{ccr} in the commutative case.

\begin{prp}
  \label{prp:glaz}
  Let $R \subseteq S$ be right coherent rings such that $S$ is
  faithfully flat as an $\Rop$-module. If $S$ is right regular, then
  $R$ is right regular.
\end{prp}

\begin{prf*}
  Let $F$ be an acyclic complex of flat $R$-modules. As $S$ is right
  regular, the acyclic complex $\tp{S}{F}$ of flat $S$-modules is pure
  acyclic by \thmref{P}; this is condition \tagFp. The character
  complex
  \begin{equation*}
    \Hom[\ZZ]{\tp{S}{F}}{\QZ} \dis \Hom{F}{\Hom[\ZZ]{S}{\QZ}}
  \end{equation*}
  is contractible. By \prpcite[1.3.49]{dcmca} the $R$-module
  $\Hom[\ZZ]{S}{\QZ}$ is faithfully injective, so $F$ is pure acyclic,
  whence $R$ is right regular by \thmref{P}.
\end{prf*}

\begin{exa}
  \label{exa:vnr}
  Let $R$ be von Neumann regular, that is, every $R$-module is
  flat. It follows that a product of flat $R$-modules is flat, whence
  $R$ is right coherent.  Flatness of every $R$-module also means that
  $R$ satisfies condition $\tagFp$, whence $R$ is right regular by
  \thmref{P}. As von Neumann regularity is a left--right symmetric
  property, $R$ is also left coherent and left regular. One could
  reach the same conclusion by recalling that every finitely generated
  left/right ideal of a von Neumann regular ring $R$ is a direct summand of $R$.
\end{exa}

A von Neumann regular ring is a special case of a right coherent ring
of finite weak global dimension; any such ring evidently satisfies
condition \tagFp\ and is thus right regular. Finkel Jones and Teply
\cite{MFJMLT82} give examples of such rings. Glaz \seccite[6.2]{ccr}
shows that the polynomial algebra in countably many variables over a
field is a coherent regular ring of infinite weak global dimension.

\begin{cor}
  \label{cor:vnr1}
  The following conditions are equivalent.
  \begin{eqc}
  \item $R$ is von Neumann regular.
  \item Every complex of finitely presented $R$-modules is
    semi-projective.
  \item Every $R$-complex is semi-flat.
  \item Every acyclic $R$-complex is pure acyclic.
  \end{eqc}
\end{cor}

\begin{prf*}
  Condition \eqclbl{iii} implies \eqclbl{iv} by \pgref{flat}(a).
  
  \proofofimp{i}{ii} Every $R$-module is flat, so every finitely
  presented $R$-module is projective. Further, $R$ is per \exaref{vnr}
  right regular, so by \thmref{P} every complex of projective
  $R$-modules is semi-projective.

  \proofofimp{ii}{iii} An $R$-complex is a filtered colimit of
  complexes of finitely presented $R$-modules, see
  \prpcite[3.3.21]{dcmca}, i.e.\ a filtered colimit of semi-flat
  $R$-complexes, see \pgref{flat}(c), and hence semi-flat by
  \prpcite[5.4.13]{dcmca}.

  \proofofimp{iv}{i} Let $M$ be an $R$-module. There is a projective
  $R$-module $P$ and an exact sequence $0 \to K \to P \to M \to 0$; as
  it is pure, $M$ is flat.
\end{prf*}

\section{Complexes of flat-cotorsion modules}
\label{sec:2}

\noindent A cotorsion pair $(\mathsf{X},\mathsf{Y})$ in the category
$\Mod$ of $R$-modules induces two cotorsion pairs $(\ac{X},\semi{Y})$
and $(\semi{X}, \ac{Y})$ in the category of $R$-complexes; this was
proved by Gillespie \prpcite[3.6]{JGl04}. Here the class $\ac{X}$
consists of acyclic complexes $X$ with cycle modules $\Cy[n]{X}$ from
$\mathsf{X}$, while $\semi{Y}$ consists of complexes $Y$ of modules
from $\mathsf{Y}$ with the property that $\Hom{X}{Y}$ is acyclic for
every complex $X$ in $\ac{X}$. Similarly $\ac{Y}$ consists of acyclic
complexes $Y$ with cycle modules $\Cy[n]{Y}$ from $\mathsf{Y}$, and
$\semi{X}$ consists of complexes $X$ of modules from $\mathsf{X}$ with
the property that $\Hom{X}{Y}$ is acyclic for every complex $Y$ from
$\ac{Y}$. If the cotorsion pair $(\mathsf{X},\mathsf{Y})$ is complete
and hereditary, then the induced cotorsion pairs are both complete;
see Yang and Liu \thmcite[3.5]{GYnZLi11}.

Let $\Prj$, $\Inj$, $\Flat$, and $\Cot$ be the classes of projective,
injective, flat, and cotorsion $R$-modules. For each of the cotorsion
pairs $(\Prj,\Mod)$ and $(\Mod,\Inj)$ only one of the induced
cotorsion pairs in the category of $R$-complexes is of interest, and
they yield the notions of semi-projective and semi-injective
complexes.  We proceed to recall the key properties of the cotorsion
pairs induced by the complete hereditary cotorsion pair $(\Flat,\Cot)$
in $\Mod$. A module in $\Flat \cap \Cot$ is
called~\emph{flat-cotorsion.}

\begin{lem}
  \label{lem:flatcot}
  For every $R$-complex $X$ there are exact sequences of $R$-complexes
  \begin{equation}
    \label{eq:flatcot}
    0 \lra C' \lra F \lra X \lra 0 \quad\text{ and }\quad
    0 \lra X \lra C \lra F' \lra 0      \,.
  \end{equation}
  Here $F$ is semi-flat and $C'$ is an acyclic complex of cotorsion
  modules, while $C$ is a complex of cotorsion modules, and $F'$ is
  acyclic and semi-flat.
  \begin{prt}
  \item An $R$-complex is semi-flat if and only if it belongs to
    $\semi{\Flat}$, and a complex in $\ac{\Flat}$ is a pure acyclic
    complex of flat $R$-modules.
  \item Every complex of cotorsion $R$-modules belongs to
    $\semi{\Cot}$ and every acyclic complex of cotorsion $R$-modules
    belongs to $\ac{\Cot}$.
  \item For an acyclic complex $F$ of flat-cotorsion $R$-modules the
    next conditions are equivalent.
    \begin{eqc}
    \item $F$ is semi-flat.
    \item $F$ is contractible.
    \item $F$ is pure acyclic.
    \end{eqc}
  \end{prt}
\end{lem}

\begin{prf*}
  An acyclic complex with flat cycle modules is a pure acyclic complex
  of flat $R$-modules, and such a complex is semi-flat by
  \pgref{flat}(a). The remaining assertions in parts (a) and (b) hold
  by \thmcite[1.3]{BCE-20} and \prpcite[1.6]{CELTWY-21}. It follows
  that the exact sequences in \eqref{flatcot} are standard
  approximation sequences associated to the induced cotorsion
  pairs. In part (c) conditions \eqclbl{i} and \eqclbl{iii} are
  equivalent by \pgref{flat}(a). Evidently, \eqclbl{ii} implies
  \eqclbl{iii}, and the converse follows from part (b). Indeed, in a
  pure acyclic complex of flat modules the cycle modules are flat and
  by (b) they are cotorsion as well.
\end{prf*}
  
\begin{prp}
  \label{prp:FC}
  The following conditions are equivalent and equivalent to conditions
  \tagsP\ and \tagsF\ from {\rm\prpref[]{P}}.
  \setlength{\thmlistleft}{3.25em} %
  \begin{rqm}
  \item[\tagFC] Every complex of flat-cotorsion $R$-modules is
    semi-flat.
  \item[\tagFCc] Every acyclic complex of flat-cotorsion $R$-modules
    is contractible.
  \item[\tagFCa] Every acyclic complex of flat-cotorsion $R$-modules
    is semi-flat.
  \end{rqm}
  \setlength{\thmlistleft}{2.5em} %
\end{prp}

\begin{prf*}
  By \lemref{flatcot}(c) conditions \tagFCc\ and \tagFCa\ are
  equivalent; \tagFC\ clearly implies \tagFCa. For the converse let
  $X$ be a complex of flat-cotorsion $R$-modules and consider the
  exact sequence $0 \to C' \to F \to X \to 0$ from \eqref{flatcot}. It
  follows that $C'$ is a complex of flat-cotorsion modules and, hence,
  semi-flat, so $X$ is semi-flat by \pgref{flat}(b). Thus, conditions
  \tagsFC\ are equivalent.

  Evidently \tagF\ implies \tagFC. For the
  converse let $X$ be a complex of flat $R$-modules and consider the
  exact sequence $0 \to X \to C \to F' \to 0$ from \eqref{flatcot}. It
  follows that $C$ is a complex of flat-cotorsion modules and, hence,
  semi-flat. Now $X$ is semi-flat by \pgref{flat}(b).
\end{prf*}

\begin{rmk}
  \label{rmk:3}
  Assume that $R$ is right coherent. For an acyclic complex $I$ of
  injective $\Rop$-modules, the character complex $\Hom[\ZZ]{I}{\QZ}$
  is an acyclic complex of flat-cotorsion $R$-modules. If it is
  contractible, then $I$ is pure acyclic and, therefore, contractible
  by \pgref{inj}(a). Thus \tagFCc\ implies \tagIc, i.e.\ conditions
  \tagsP, \tagsF, \tagsI, and \tagsFC\ are per \rmkref{0} and
  \prpref{FC} equivalent.
\end{rmk}

We can now add a condition to the characterization of von Neumann
regular rings. Over such a ring every cotorsion module is injective,
so \eqclbl{v} below can be seen as the counterpart to $\tagIc$ in the
characterization of regular rings.

\begin{cor}
  \label{cor:vnr2}
  The next condition is equivalent to conditions
  \eqclbl{i}--\eqclbl{iv} from {\rm\corref[]{vnr1}}.
  \begin{eqc}
    \setcounter{eqc}{4}
  \item Every acyclic complex of cotorsion $R$-modules is
    contractible.
  \end{eqc}
\end{cor}

\begin{prf*}
  Recall the equivalent conditions \eqclbl{i}--\,\eqclbl{iv} from
  \corref{vnr1}.

  \proofofimp{i}{v} Every $R$-module is flat, and $R$ is right
  coherent and right regular by \exaref{vnr}, so every acyclic complex
  of cotorsion $R$-modules is a complex of flat-cotorsion modules and
  hence contractible by \thmref{P} and \prpref{FC}.

  \proofofimp{v}{i} The augmented injective resolution of a cotorsion
  $R$-module is contractible, so every cotorsion $R$-module is
  injective and every $R$-module is flat.
\end{prf*}

\section{Complexes of fp-injective modules}
\label{sec:3}

\noindent
Recall that an $\Rop$-module $E$ is fp-injective if
$\Ext[\Rop]{1}{F}{E}=0$ holds for every finitely presented
$\Rop$-module $F$. The fp-injective $\Rop$-modules constitute the
right-hand class of a cotorsion pair; the modules in the left-hand
class are known as fp-projective. If $R$ is right coherent, then this
cotorsion pair, $(\FPPrj,\FPInj)$, is complete and hereditary and in
many ways dual to $(\Flat,\Cot)$. In the next lemma we collect the key
properties of the induced cotorsion pairs in the category of
$R$-complexes. We refer to complexes in $\semi{\FPInj}$ as
\emph{semi-fp-injective} and to modules in $\FPPrj \cap \FPInj$ as
\emph{fp-pro-injective.}

\begin{lem}
  \label{lem:fp}
  Let $R$ be right coherent. For every $\Rop$-complex $X$ there are
  exact sequences of $\Rop$-complexes
  \begin{equation}
    \label{eq:fp}
    0 \lra X \lra E \lra P' \lra 0   \quad\text{ and }\quad
    0 \lra E' \lra P \lra X \lra 0 \,.       
  \end{equation}
  Here $E$ is semi-fp-injective and $P'$ is an acyclic complex of
  fp-projective modules, while $P$ is a complex of fp-projective
  modules and $E'$ is acyclic and semi-fp-injective.
  \begin{prt}
  \item Every complex of fp-projective $\Rop$-modules belongs to
    $\semi{\FPPrj}$ and every acyclic complex of fp-projective
    $\Rop$-modules is in $\ac{\FPPrj}$.
  \item An acyclic complex of fp-injective $\Rop$-modules is pure
    acyclic if and only if it is semi-fp-injective.
  \item If in an exact sequence $0 \to E' \to E \to E'' \to 0$ of
    complexes of fp-injective $\Rop$-modules two of the complexes are
    semi-fp-injective, then so is the third.
  \item For an acyclic complex $E$ of fp-pro-injective $\Rop$-modules
    the next conditions are equivalent.
    \begin{eqc}
    \item $E$ is semi-fp-injective.
    \item $E$ is contractible.
    \item $E$ is pure acyclic.
    \end{eqc}
  \end{prt}
\end{lem}

\begin{prf*}
  It follows from parts (a) and (b) that the exact sequences in
  \eqref{fp} are standard approximation sequences associated to the
  induced cotorsion pairs.
  
  The assertions in part (a) were proved by \v{S}aroch and
  \v{S}tov{\'{\i}}{\v{c}}ek \exacite[4.3]{JSrJSt20}.

  \proofoftag{b} Let $E$ be an acyclic complex of fp-injective
  $\Rop$-modules. If $E$ is semi-fp-injective, then each cycle module
  $\Cy[n]{E}$ is fp-injective, see Gillespie \corcite[3.13(5)]{JGl04}
  and, therefore, $E$ is pure acyclic. Conversely, if $E$ is pure
  acyclic, then the cycle modules $\Cy[n]{E}$ are fp-injective and $E$
  is semi-fp-injective by \lemcite[3.10]{JGl04}.

  \proofoftag{c} Let $P$ be an acyclic complex of fp-projective
  $\Rop$-modules. There is an induced exact sequence
  \begin{equation*}
    0 \lra \Hom[\Rop]{P}{E'} \lra \Hom[\Rop]{P}{E} \lra \Hom[\Rop]{P}{E''} \lra 0
  \end{equation*}
  and if two of these complexes are acyclic, then so is the third.

  \proofoftag{d} Let $E$ be an acyclic complex of fp-pro-injective
  $\Rop$-modules. Conditions \eqclbl{i} and \eqclbl{iii} are
  equivalent by \prtlbl{b}, and every contractible complex is pure
  acyclic. If $E$ is pure acyclic, then the cycle modules $\Cy[n]{E}$
  are fp-injective, and by part (a) they are fp-projective as well, so
  $E$ is contractible.
\end{prf*}

\pagebreak
\begin{thm}
  \label{thm:FP}
  Let $R$ be right coherent. The following conditions are equivalent,
  and equivalent to conditions \tagsP, \tagsF, \tagsI, and \tagsFC\
  from {\rm\prpref[]{P}}, {\rm\prpref[]{I}}, and {\rm\prpref[]{FC}}.
  \setlength{\thmlistleft}{3.75em}%
  \begin{rqm}
  \item[\tagFPI] Every complex of fp-injective $\Rop$-modules is
    semi-fp-injective.
  \item[\tagFPIp] Every acyclic complex of fp-injective $\Rop$-modules
    is pure acyclic.
  \item[\tagFPIa] Every acyclic complex of fp-injective $\Rop$-modules
    is semi-fp-injective.

  \item[\tagFPPI] Every complex of fp-pro-injective $\Rop$-modules is
    semi-fp-injective.
  \item[\tagFPPIc] Every acyclic complex of fp-pro-injective
    $\Rop$-modules is contractible.
  \item[\tagFPPIa] Every acyclic complex of fp-pro-injective
    $\Rop$-modules is semi-fp-injective.
  \end{rqm}
  \setlength{\thmlistleft}{2.5em} %
\end{thm}

\begin{prf*}
  First we argue that the six conditions \tagsFPI\ and \tagsFPPI\ are
  equivalent. By \lemref{fp}(b,d) one has
  \begin{equation*}
    \tagFPI \implies \tagFPIp \implies \tagFPIa \qquad\text{and}\qquad
    \tagFPPI \implies \tagFPPIc \implies \tagFPPIa \,,
  \end{equation*}
  which leaves two implications to prove.

  \tagFPIa$\implies$\tagFPPI: Let $X$ be a complex of fp-pro-injective
  $\Rop$-modules and consider the exact sequence
  $0 \to X \to E \to P' \to 0$ from \eqref{fp}. It follows that $P'$
  is a complex of fp-injective $\Rop$-modules and, therefore,
  semi-fp-injective, so $X$ is semi-fp-injective by \lemref{fp}(c).

  \tagFPPIa$\implies$\tagFPI: Let $X$ be a complex of fp-injective
  $\Rop$-modules and consider the exact sequence
  $0 \to X \to E \to P' \to 0$ from \eqref{fp}. It follows that $P'$
  is an acyclic complex of fp-pro-injective $\Rop$-modules and,
  therefore, semi-fp-injective, so $X$ is semi-fp-injective by
  \lemref{fp}(c).

  To finish, it suffices by \rmkref{0} and \prpref{FC} to prove two
  implications.

  \tagFPIp$\implies$\tagIc: Let $I$ be an acyclic complex of injective
  $\Rop$-modules. It is pure acyclic by \tagFPIp\ and hence
  contractible by \pgref{inj}(a).

  \tagFCc$\implies$\tagFPIp: Let $E$ be an acyclic complex of
  fp-injective $\Rop$-modules. The character complex
  $\Hom[\ZZ]{E}{\QZ}$ is an acyclic complex of flat-cotorsion
  $R$-modules and hence contractible, so $E$ is pure acyclic.
\end{prf*}

\begin{cor}
  \label{cor:conditions}
  Let $R$ be right coherent; it is right regular if and only if
  one/all of the nineteen conditions \tagsP, \tagsF, \tagsI, \tagsFC,
  \tagsFPI, and \tagsFPPI\ from {\rm\prpref[]{P}}, {\rm\prpref[]{I}},
  {\rm\prpref[]{FC}}, and {\rm\thmref[]{FP}} hold.
\end{cor}

\begin{prf*}
  Combine \thmref[Theorems~]{P} and \thmref[]{FP}.
\end{prf*}

\begin{rmk}
  \label{rmk:2}
  Gillespie and Iacob \cite[Thms.~4.3--4.6]{JGlAIc} characterize right
  coherent right regular rings in terms of the equivalence of
  conditions \tagP, \tagPc, \tagF, \tagFp, \tagI, \tagIc, \tagFPI, and
  \tagFPIp\ from {\rm\prpref[]{P}}, {\rm\prpref[]{I}}, and
  {\rm\thmref[]{FP}}, among other conditions not considered here.
\end{rmk}

Finally we can add conditions to the characterization of von Neumann
regular rings. Over such a ring every fp-projective module is
projective, so \eqclbl{vii} below can be seen as the counterpart to
$\tagPc$ in the characterization of regular rings.

\begin{cor}
  The following conditions are equivalent and equivalent to conditions
  \eqclbl{i}--\eqclbl{v} from {\rm\corref[]{vnr1}} and
  {\rm\corref[]{vnr2}}.
  \begin{eqc}
    \setcounter{eqc}{5}
  \item Every $R$-complex is semi-fp-injective.
  \item Every acyclic complex of fp-projective $R$-modules is
    contractible.
  \end{eqc}
\end{cor}

\begin{prf*}
  Recall the equivalent conditions \eqclbl{i}--\,\eqclbl{v} from
  \corref[Corollaries~]{vnr1} and \corref[]{vnr2}.

  \proofofimp{i}{vii} A finitely presented $R$-module is by
  \corref{vnr1} projective, so every $R$-module is fp-injective. Thus,
  an fp-projective module is an fp-pro-injective $R$-module. Further,
  $R$ is by \exaref{vnr} left coherent and left regular, so every
  acyclic complex of fp-pro-injective $R$-modules is contractible by
  \corref{conditions}.
  
  \proofofimp{vii}{vi} Let $M$ be an $R$-complex. A complex $P$ in
  $\ac{\FPPrj[R]}$ is, in particular, an acyclic complex of
  fp-projective $R$-modules and, therefore, contractible. Thus
  $\Hom{P}{M}$ is acyclic, whence $M$ is semi-fp-injective.
  
  \proofofimp{vi}{i} As every $R$-module is fp-injective, every
  finitely presented $R$-module is projective. Since a module is a
  filtered colimit of its finitely presented submodules, it follows
  that every $R$-module is flat, so $R$ is von Neumann regular.
\end{prf*}

\def\soft#1{\leavevmode\setbox0=\hbox{h}\dimen7=\ht0\advance \dimen7
  by-1ex\relax\if t#1\relax\rlap{\raise.6\dimen7
  \hbox{\kern.3ex\char'47}}#1\relax\else\if T#1\relax
  \rlap{\raise.5\dimen7\hbox{\kern1.3ex\char'47}}#1\relax \else\if
  d#1\relax\rlap{\raise.5\dimen7\hbox{\kern.9ex \char'47}}#1\relax\else\if
  D#1\relax\rlap{\raise.5\dimen7 \hbox{\kern1.4ex\char'47}}#1\relax\else\if
  l#1\relax \rlap{\raise.5\dimen7\hbox{\kern.4ex\char'47}}#1\relax \else\if
  L#1\relax\rlap{\raise.5\dimen7\hbox{\kern.7ex
  \char'47}}#1\relax\else\message{accent \string\soft \space #1 not
  defined!}#1\relax\fi\fi\fi\fi\fi\fi}
  \providecommand{\MR}[1]{\mbox{\href{http://www.ams.org/mathscinet-getitem?mr=#1}{#1}}}
  \renewcommand{\MR}[1]{\mbox{\href{http://www.ams.org/mathscinet-getitem?mr=#1}{#1}}}
  \providecommand{\arxiv}[2][AC]{\mbox{\href{http://arxiv.org/abs/#2}{\sf
  arXiv:#2 [math.#1]}}} \def\cprime{$'$}
\providecommand{\bysame}{\leavevmode\hbox to3em{\hrulefill}\thinspace}
\providecommand{\MR}{\relax\ifhmode\unskip\space\fi MR }
\providecommand{\MRhref}[2]{%
  \href{http://www.ams.org/mathscinet-getitem?mr=#1}{#2}
}
\providecommand{\href}[2]{#2}

\end{document}